\documentclass{article}
\usepackage{amsmath}
\usepackage{amssymb}
 \newtheorem{thm}{Theorem}[section]
 \newtheorem{cor}[thm]{Corollary}
 
 \newtheorem{prop}[thm]{Proposition}

 \newtheorem{rem}[thm]{Remark}
 \newtheorem{nota}[thm]{Notation}
\newtheorem{quest}[thm]{Question}
 \newtheorem{defn}[thm]{Definition}
 \newtheorem{ex}[thm]{Example}
 \numberwithin{equation}{section}
\newcommand{\spec}{{\rm Spec}}
\newcommand{\de}{{\rm depth}}
\newcommand{\aCM}{{\rm {\bf aCM}}}
\newcommand{\supp}{{\rm Supp}}
\newcommand{\hgt}{{\rm ht}}
\newcommand{\smin}{{\rm Min}}
\newcommand{\ass}{{\rm Ass}}
\newcommand{\spm}{{\rm Max}}
    \begin{document}
  \title{More properties  of almost Cohen-Macaulay rings}
  \author{Cristodor Ionescu\\[2mm]
  Institute of Mathematics \textit{Simion Stoilow}\\ of the Romanian Academy\\
  P.O. Box 1-764\\ RO 014700 Bucharest\\  Romania
  \\
  email: cristodor.ionescu@imar.ro}
  \date{}
  \maketitle

  \begin{abstract}
Some interesting properties of almost Cohen-Macaulay rings are investigated and a Serre  type property connected with this class of rings is studied. 
\end{abstract}
\vspace{5mm}
\section{Introduction}
A flaw in the chapter dedicated to Cohen-Macaulay rings in the first edition of \cite{Mat} was corrected in the second edition. This led to the study  of the so-called almost Cohen Macaulay rings, 
first by Y. Han \cite{Han} and later by M.-C. Kang \cite{K1}, \cite{K2}. Since the first of these papers is written in Chinese, the others two are the main reference for the subject. 

\begin{rem}\label{fordef}
Let $A$ be a commutative noetherian ring, $P\in\spec(A)$  and $M\neq 0$ a finitely generated $A$-module. Then $\de_P(M)\leq\de_{PA_P}M_P.$ 
\end{rem}

\begin{defn}\label{maindef} {\rm (cf. \cite{Han}, \cite{K1})}
Let $A$ be a commutative noetherian ring. A finitely generated $A$-module $M\neq 0$ is called almost Cohen-Macaulay if
$\de_PM=\de_{PA_P}M_P$, for any $P\in\supp(M).$ $A$ is called  an almost Cohen-Macaulay ring if it is an almost Cohen-Macaulay  $A$-module, that is if for any $P\in\spec(A),$ $\de_PA=\de_{PA_P}A_P.$  
\end{defn}
\par\noindent Several properties of almost Cohen-Macaulay rings are proved in \cite{K1} and several interesting examples are given in \cite{K2}. In the following we are trying to complete the results
in \cite{K1} and to introduce a Serre-type condition, that we call $(C_k),$ for any $ k\in\mathbb{N},$ condition that is to be to almost Cohen-Macaulay rings what the classical Serre condition $(S_k)$ is to 
Cohen-Macaulay rings. 
\section{Properties of almost Cohen-Macaulay rings}

All rings considered will be commutative and with unit. We start by reminding some basic properties of almost Cohen-Macaulay rings.

\begin{rem}\label{obs1} Let $A$ be a noetherian ring. Then:
\par\noindent a) $A$ is almost Cohen-Macaulay iff $\hgt(P)\leq 1+\de_PA, \forall\ P\in\spec(A)$ {\rm (\cite{K1}, 1.5)};
\par\noindent b) $A$ is almost Cohen-Macaulay iff $A_P$ is almost Cohen-Macaulay for any $P\in\spec(A)$ iff $A_Q$ is almost Cohen-Macaulay for any $Q\in\spm(A)$ iff $\hgt(Q)\leq 1+\de A_Q$ for any $Q\in\spm(A)$ {\rm (\cite{K1}, 2.6)};
\par\noindent c) If $A$ is local, it follows from b) that $A$ is almost Cohen-Macaulay if and only if $\dim(A)\leq 1+\de(A).$
\end{rem}

\par\noindent Our first result is a stronger formulation of \cite{K1}, 2.10 and deals with the behaviour of almost Cohen-Macaulay rings with respect to flat morphisms.

\begin{prop}\label{flat}
Let $u:(A,m)\to(B,n)$ be a local flat morphism of noetherian local rings. 
\par\noindent a) If $B$ is almost Cohen-Macaulay, then $A$ and $B/mB$ are almost Cohen-Macaulay.
\par\noindent b)If $A$ and $B/mB$ are almost Cohen-Macaulay and one of them is Cohen-Macaulay, then $B$ is almost Cohen-Macaulay.
\end{prop}

\par\noindent\textit{Proof:}
a) We have $$\dim(A)=\dim(B)-\dim(B/mB)\leq 1+\de B-\dim(B/mB)\leq$$
$$\leq 1+\de B-\de(B/mB)=1+\de A.$$
We have also  
$$\dim(B/mB)-\de(B/mB)=(\dim(B)-\de B)-(\dim(A)-\de A)\leq$$
$$\leq 1-(\dim(A)-\de A)\leq 1.$$
\par\noindent b) Since $u$ is flat we have $$\dim(B)=\dim(A)+\dim(B/mB)\leq 1+\de(A)+\de (B/mB)=$$
$$=1+\de(B).$$

\begin{quest}\label{prob}

We don't know of any example of a   local flat morphism of noetherian local rings $u:(A,m)\to(B,n)$ such that  $A$ and $B/mB$ are almost Cohen-Macaulay 
and $B$ is not almost Cohen-Macaulay.

\end{quest}
\begin{cor}\label{compl}
Let $A$ be a noetherian local ring, $I\neq A$ be an ideal contained in the Jacobson radical of $A$ and  $\hat{A}$ the completion of $A$ in the $I$-adic topology. 
Then $A$ is almost Cohen-Macaulay if and only if $\hat{A}$  is almost Cohen-Macaulay.
\end{cor}
\par\noindent\textit{Proof:}
Since $I$ is contained in the Jacobson radical of $A,$ the canonical morphism $A\to\hat{A}$ is faithfully flat and $\spm(A)\cong\spm(\hat{A}).$ Moreover, if $m\in\spm(A)$ and $\hat{m}$ is the 
corresponding maximal ideal of $\hat{A},$ the closed fiber of the morphism $A_m\to\hat{A}_{\hat{m}}$ is a field. Now apply \ref{flat}.

\begin{cor}\label{series} {\rm (see \cite{K1}, 1.6)}
Let $A$ be a noetherian ring and $n\in\mathbb{N}.$ Then $A$ is almost Cohen-Macaulay if and only if $A[[X_1,\ldots,X_n]]$ is almost Cohen-Maculay.
\end{cor}

\par\noindent\textit{Proof:} 
Suppose  that $A$ is almost Cohen-Macaulay. We may clearly assume that $A$ is local and $n=1.$ By \cite{K1}, 1.3 we get that $A[X]_{(X)}$ is almost Cohen-Macaulay. Now apply \ref{compl}.
The converse is clear.  

\par For the next corollary we need some notations.
\begin{nota}\label{notatie} 
If \textbf{P} is a property of noetherian local rings, we denote by $\textbf{P}(A):=\{Q\in\spec(A)\mid A_Q\ \text{has the property}\ \textbf{P}\}$ and by $\textbf{NP}(A):=\{Q\in\spec(A)\mid A_Q\ \text{has not the property}\ \textbf{P}\}=\spec(A)\setminus\textbf{P}(A).$
\end{nota}

\begin{defn}\label{acm}
Let $A$ be a noetherian ring. According to \ref{notatie}, the set 
$$\aCM(A):=\{P\in\spec(A)\ \vert\ A_P\ \text{is almost Cohen-Macaulay}\}$$
is called the almost Cohen-Macaulay locus of $A.$
\end{defn}

\begin{cor}\label{spec}
Let $u:A\to B$ be a morphism of noetherian local rings and $\varphi:\spec(B)\to \spec(A)$ the induced morphism on the spectra. If the fibers of $u$  are Cohen-Macaulay, 
then $\varphi^{-1}(\aCM(A))=\aCM(B).$
\end{cor}

\par\noindent\textit{Proof:}
Obvious from \ref{flat}.

In Cohen-Macaulay rings chains of prime ideals behave very well, in the sense that Cohen-Macaulay rings are universally catenary(see \cite{Mat}). This is no more the case for almost Cohen-Macaulay rings.  

\begin{ex}\label{uncat}
There exists a local  almost Cohen-Macaulay ring which is not  catenary.
\end{ex}

\par\noindent\textit{Proof:}
Indeed, by \cite{K1}, Ex. 2, any noetherian normal integral domain of dimension 3 is almost Cohen-Macaulay. In \cite{Og} such a ring which is not catenary is constructed.

The next result shows that some of the formal fibres of almost Cohen-Macaulay rings are almost Cohen-Macaulay. A stronger fact will be proved in \ref{ffiber}.

\begin{prop}\label{brh}
Let $A$ be a noetherian local almost Cohen-Macaulay ring, $P\in\spec(A), Q\in\ass(\hat{A}/P\hat{A}).$ Then $\hat{A}_Q/P\hat{A}_Q$ is almost Cohen-Macaulay.
\end{prop}

\par\noindent\textit{Proof:}
We have 
$$\dim(\hat{A}_Q/P\hat{A}_Q)=\dim\hat{A}_Q-\dim A_P\leq\de\hat{A}_Q+1-\dim A_P\leq$$
$$\leq\de\hat{A}_P+1-\dim A_P=\de(\hat{A}_Q/P\hat{A}_Q)+1.$$

The following result shows that the almost Cohen-Macaulay property is preserved by tensor products and finite field extensions.

\begin{prop}\label{tensor}
Let $k$ be a field, $A$ and $B$ be two $k-$algebras such that $A\otimes_k B$ is a noetherian ring. If $A$ and $B$ are almost Cohen-Macaulay, then $A\otimes_k B$ is almost Cohen-Macaulay.  
\end{prop}

\par\noindent\textit{Proof:} Let $P\in\spec(A).$ We have a flat morphism $B\to B\otimes_k k(P)$ and let $Q\in\spec(B).$ Set $T:=A/P\otimes_k B/Q=A\otimes_k B/(P\otimes_k B+A\otimes_k Q).$  
Then $k(P)\otimes_k k(Q)$ is a ring of fractions of $T,$ hence noetherian by assumption. By \cite{TY}, 1.5,\ it follows that  $k(P)\otimes_k k(Q)$ is locally a complete intersection.
Let now $Q\in\spec(B)$ and $P=Q\cap A.$ By the above  the flat local morphism $A_P\to (B\otimes_kk(P))_Q$ has a complete intersection closed fiber, hence the ring $(B\otimes_k k(P))_Q$ is almost 
Cohen-Macaulay by \ref{flat}. Now consider the flat morphism $A\to A\otimes_k B,$ let $Q\in\spec(A\otimes_k B)$ and $P=Q\cap A.$ 
Then the flat local morphism $A_P\to (A\otimes_k B)_Q$ has a complete intersection closed fiber, whence $(A\otimes_k B)_Q$ is almost Cohen-Macaulay.

\begin{cor}\label{geomacm}
Let $k$ be a field, $A$ a noetherian $k-$algebra which is almost Cohen-Macaulay and $L$ a finite field extension of $k.$ Then $A\otimes_k L$ is almost Cohen-Macaulay.
\end{cor}

As for the Cohen-Macaulay property, the formal fibres of factorizations of almost Cohen-Macaulay rings are almost Cohen-Macaulay.

\begin{prop}\label{ffiber}
Let $B$ be a local almost Cohen-Macaulay ring, $I$ an ideal of $B$ and $A=B/I.$ Then the formal fibers of $A$ are almost Cohen-Macaulay. 
\end{prop}

\par\noindent\textit{Proof:}
We have $\hat{A}=\hat{B}\otimes_B A=\hat{B}/I\hat{B},$ hence the formal fibers of $A$ are exactly
the formal fibers of $B$ in the prime ideals of $B$ containing $I.$ Let $P$ be such a
prime ideal, let $S = B\setminus P$ and let $C := S^{-1}(\hat{B} /I\hat{B}).$ Let also $Q\in\spec(C).$ There exists $Q'\in\spec(\hat{B})$ such that $Q=Q'C$ and $Q'\cap B=P.$ Thus we have a local flat
morphism $B_Q\to \hat{B}_{Q'}.$ But $B$ is almost Cohen-Macaulay, hence $\hat{B}_{Q'}$ and consequently $C_Q\cong\hat{B}_{Q'}/P\hat{B}_{Q'}$ are almost Cohen-Macaulay, by \ref{flat}.

\section{The property $(C_n)$} 

Recall that given a natural number $n,$ a noetherian ring $A$ is said to have Serre property $(S_n)$ if $\de(A_P)\geq\min(\hgt{P},n)$ for any prime ideal $P\in\spec(A).$  
Moreover, $A$ is Cohen-Macaulay if and only if $A$ has the property $(S_n)$ for any $n\in\mathbb{N}$ (see \cite{Mat}, (17.I)). We will try to characterize almost Cohen-Macaulay rings in a similar way.

\begin{defn}\label{tn}
Let $n\in\mathbb{N}$ be a natural number. We say that a noetherian ring $A$ has the property $(C_n)$ if $\de(A_P)\geq\min(\hgt{P},n)-1, \forall\ P\in\spec(A).$ 
\end{defn}

\begin{rem}\label{tnsn}
a) It is clear that $(C_n)\Rightarrow (C_{n-1})$ and that $(S_n)\Rightarrow (C_n), \forall n\in\mathbb{N}.$
\par\noindent b) It is also clear that if $A$ has $(C_n),$ then $A_P$ has $(C_n), \forall P\in\spec(A).$
\end{rem}

\begin{thm}\label{acmtn}
A noetherian ring $A$ is almost Cohen-Macaulay if and only if $A$ has the property $(C_n)$ for every $n\in\mathbb{N}.$
\end{thm}

\par\noindent\textit{Proof:}
Assume that $A$ is almost Cohen-Macaulay and let $P\in\spec(A).$ Then $A_P$ is almost Cohen-Macaulay, hence $\de(A_P)\geq\hgt(P)-1.$ If $n\geq\hgt(P),$ then $\min(\hgt(P),n)=\hgt(P),$ hence $\de(A_P)\geq\min(n,\hgt(P))-1.$ If $n<\hgt(P),$ then $\min(n,\hgt(P))=n,$ so that $\de(A_P)\geq\hgt(P)-1>n-1=\min(\hgt(P),n)-1.$
\par\noindent For the converse, let $P\in\spec(A), \hgt(P)=l.$ Then 
\[
\de(A_P)\geq\min(l,\hgt(P))-1=\hgt(P)-1.
\]

\begin{prop}\label{loctk}
Let $k\in\mathbb{N}.$ A noetherian ring $A$ has the property $(C_k)$ if and only if $A_P$ is almost Cohen-Macaulay for any $P\in\spec(A)$ with $\de(A_P)\leq k-2.$
\end{prop}

\par\noindent\textit{Proof:}
Let $P\in\spec(A)$ such that $\min(k,\hgt(P))-1\leq\de(A_P)\leq k-2.$ If $\hgt(P)\leq k,$ then $\de(A_P)\geq\hgt(P)-1.$ And if $\hgt(P)>k,$ then it follows that $k-2>\de(A_P)\geq k-1.$ Contradiction!
\par\noindent Conversely, let $P\in\spec(A).$ If $\de(A_P)\leq k-2,$ then $A_P$ is almost Cohen-Macaulay, hence $\hgt(P)-1\leq\de(A_P)\leq k-2.$ Thus $\min(\hgt(P),k)=\hgt(P),$ whence $\de(A_P)\geq\min(k,\hgt(P).$ If $k-2<\de(A_P),$ then $\hgt(P)>k-2,$ hence $\de(A_P)\geq\min(k,\hgt(P))-1.$

\begin{prop}\label{regel}
Let $A$ be a noetherian ring, $k\in\mathbb{N}$ and $x\in A$ a non zero divisor. If $A/xA$ has the property $(C_k)$, then $A$ has the property $(C_k).$ 
\end{prop}

\par\noindent\textit{Proof:}
Let $Q\in\spec(A)$ such that $\de(A_Q)=n\le k-2.$ If $x\in Q,$ then $\de(A/xA)_Q=n-1\leq k-3.$ Then $\hgt(Q/xA)\leq n-1+1=n,$ hence $\hgt(Q)\leq n+1=\de A_Q+1.$  If $x\notin Q,$ let $P\in\smin(Q+xA).$ Then $(P+xA)A_Q$ is $QA_Q$-primary and $\de(A_P)\leq\de(A_Q)+1=n+1.$ Then $\de(A/xA)_Q=n-1,$ hence $\hgt(P/xA)\leq n.$ It follows that $\hgt(P)\leq n+1=\de(A_P)+1.$

\begin{defn}\label{nagcr} 
We say that a property \textbf{P} of noetherian local rings satisfies Nagata's Criterion (NC) if the following holds:
if $A$ is a noetherian ring such for every $P\in \textbf{P}(A),$ the set $\textbf{P}(A/P)$ contains a non-empty open set of $\spec(A/P)$, then $\textbf{P}(A)$ is open in $\spec(A).$ 
\end{defn}

\par\noindent An interesting study of Nagata Criterion is performed in \cite{Mas}.

\begin{thm}\label{Tknagata}
Let $k\in\mathbb{N}.$ The property $(C_k)$ satisfies (NC).
\end{thm}

\par\noindent\textit{Proof:} 
Let $Q\in C_k(A).$ Then $\de(A_Q)\geq\min(k,\hgt(Q))-1.$
\par\noindent Case a): $\hgt(Q)\leq k.$ Then $\min(k,\hgt(Q))=\hgt(Q),$ hence $\de(A_Q)+1\geq\hgt(Q)$ and $A_Q$ is almost Cohen-Macaulay. Let $f\in A\setminus Q$ such that 
$$\dim(A_P)=\dim(A_Q)+\dim(A_P/QA_P)$$
 and 
 $$\de(A_P)=\de(A_Q)+\de(A_P/QA_P)$$
for any $P\in D(f)\cap V(Q)\cap NT_k(A).$ Then $\de(A_P)\ngeqq\min(k,\hgt(P))-1.$
\par\noindent Case a1): $\hgt(P)\leq k.$ Then $\min(k,\hgt(P))=\hgt(P),$ hence $\de(A_P)+1<\hgt(P).$ Then $$\de(A_P/QA_P)+1=\de(A_P)-\de(A_Q)+1<$$
$$<\hgt(P)-\de(A_Q)\leq\hgt(P)-\hgt(Q)+1.$$
Then $\de(A_P/QA_P)<\dim(A_P/QA_P)=\dim(A_P)-\dim(A_Q)$ and it follows that $A_P/QA_P$ is not $(C_k).$
\par\noindent Case a2): $\hgt(P)> k.$ Then $\min(k,\hgt(P))=k,$ hence $\de(A_P)<k-1.$ It follows that 
$$\de(A_P/QA_P)=\de(A_P)-\de(A_Q)<$$
$$<k-1+1-\hgt(Q)=k-\hgt(Q).$$
This implies that $A_P/QA_P$ is not $(C_k).$  
\par\noindent Case b): $\hgt(Q)>k.$ Then $\min(k,\hgt(Q))=k$ and $\de(A_Q)+1\geq k.$ Since $\hgt(P)>k,$ it follows that $\min(k,\hgt(P))=k$ and $\de(A_P)+1<k.$ Let $x_1,\ldots,x_r$ be an 
$A_Q$-regular sequence. Then there exists $f\in A\setminus Q$ such that $x_1,\ldots,x_r$ is $A_f$-regular. If $P\in D(f)\cap V(Q),$ it follows that $A_P$ is $(C_k).$

\begin{cor}\label{acmnagata}
The property almost Cohen-Macaulay satisfies (NC).
\end{cor}

\begin{thm}\label{open}
Let $A$ be a quasi-excellent ring and $k\in\mathbb{N}.$  Then $C_k(A)$ and $\aCM(A)$ are open in the Zariski topology of $\spec(A).$
\end{thm}

\par\noindent\textit{Proof:}
Let $P\in\spec(A).$ Then $\aCM(A/P)$ and $C_k(A/P)$ contain the non-empty open set ${\bf Reg}(A/P)=\{ P\in\spec(A)\ \vert\ A_P\ \text{is regular}\ \}.$ Now apply \ref{Tknagata} and \ref{acmnagata}.

\begin{cor}\label{opencom}
Let $A$ be a complete semilocal ring and $k\in\mathbb{N}.$  Then $C_k(A)$ and $aCM(A)$ are open in the Zariski topology of $\spec(A).$
\end{cor}

\begin{cor}\label{cmfib}
 Let $A$ be a noetherian local ring with Cohen-Macaulay formal fibers. Then $\aCM(A)$ is open.
\end{cor}

\par\noindent\textit{Proof:}
 Follows from \ref{opencom} and \ref{spec}.

\begin{prop}\label{flattk}
Let $u:A\to B$ be a flat morphism of noetherian rings and $k\in\mathbb{N}.$ If $B$ has $(C_k),$ then $A$ has $(C_k).$
\end{prop}

\par\noindent\textit{Proof:}
We may assume that $A$ and $B$ are local rings and that $u$ is local. Let $P\in\spec(A)$ and $Q\in\smin(PB).$ Then $\dim(B_Q/PB_Q)=0,$ hence 
$$\de(A_P)=\de(B_Q)\geq\min(k,\dim(B_Q))-1=$$
$$=\min(k,\dim(A_P))-1.$$

\begin{prop}\label{flattk2}
Let $u:A\to B$ be a flat morphism of noetherian rings and $k\in\mathbb{N}.$ 
\par\noindent a) If $A$ has $(C_k)$ and all the fibers of $u$ have $(S_k),$ then $B$ has $(C_k).$
\par\noindent b) If $A$ has $(S_k)$ and all the fibers of $u$ have $(C_k),$ then $B$ has $(C_k).$
\end{prop}

\par\noindent\textit{Proof:}
a) Let $Q\in\spec(B), P=Q\cap A.$ Then by flatness we have
$$\dim(B_Q)=\dim(A_P)+\dim(B_Q/PB_Q),$$  $$\de(B_Q)=\de(A_P)+\de(B_Q/PB_Q).$$
By assumption we have 
$$\de(A_P)\geq\min(k,\hgt(P))-1,$$
$$\de(B_Q/PB_Q)\geq\min(k,\dim(B_Q/PB_Q).$$
Hence we have 
$$\de(B_Q)=\de(A_P)+\de(B_Q/PB_Q)\geq$$
$$\geq\min(k,\hgt(P))-1+\min(k,\dim(B_Q/PB_Q))=\min(k,\hgt(B_Q))-1.$$
\par\noindent b) The proof is the same.

As a corollary we get a new proof of a previous result.

\begin{cor}\label{flatacm}
Let $u:A\to B$ be a flat morphism of noetherian rings. 
\par\noindent a) If $B$ is almost Cohen-Macaulay, then $A$ is almost Cohen-Macaulay.
\par\noindent b) If $A$ is almost Cohen-Macaulay and the fibers of $u$ are Cohen-Macaulay, then $B$ is almost Cohen-Macaulay.
\end{cor}

\begin{ex}\label{exemplu}
Let $k$ be a field and let $X_0,X_1,X_2,Y_1,Y_2$ be indeterminates. Set $B=k[[X_0,X_1,X_2]]/(X_0)\cap(X_0,X_1)^2\cap(X_0,X_1,X_2)^3$ and $A:=B[[Y_1,Y_2]].$ It is easy to see that $A$ is a noetherian local ring with
$\dim(A)=5, \de(A)=2.$ It is also not difficult to see that $A$ has the property $(C_3)$ and not the property $(C_4).$ Similar other examples can easily be constructed.
\end{ex}

\begin{ex}\label{altex}
 Let $k$ be a field, $X,Y$ indeterminates and consider the ring $A=k[[X,Y]]/(X^2,XY).$ Then $A$ has $(C_2)$ and not $(S_2).$
\end{ex}



  \end{document}